\documentclass{amsart}

\newtheorem{theorem}{Theorem}[section]

\newtheorem{corollary}[theorem]{Corollary}

\theoremstyle{definition}

\theoremstyle{remark}
\newtheorem{remark}[theorem]{Remark}

\numberwithin{equation}{section}

\begin{document}

\title{Nonsingular $\alpha$-rigid maps: \\  Short proof}

\author{Oleg\,N.~ Ageev
 }
\address{Department of Mathematics, Moscow State Technical
University, 2nd Baumanscaya St. 5, 105005 Moscow, Russia }
\email{ageev@mx.bmstu.ru}
\curraddr{Max Planck Institute of Mathematics,
P.O.Box 7280, D-53072 Bonn, Germany}
\email{ageev@mpim-bonn.mpg.de}
\thanks{The author was supported in part by
Max Planck Institute of Mathematics, Bonn and the Programme of
Support of Leading Scientific Schools of the RF (grant no.
NSh-3038.2008.1)}

\subjclass[2000]{Primary 37A, 28D05;
 Secondary 47A05, 47A35, 47D03, 47D07}

\date{May 03. 2009.}


\keywords{$\alpha$-rigid transformation, Lebesgue spectrum, Markov
operators}

\begin{abstract}
It is shown that for every $\alpha$, where $\alpha\in [0, 1/2]$,
there exists an $\alpha$-rigid transformation whose spectrum has
Lebesgue component. This answers the question posed by Klemes and
Reinhold in [7]. We apply a certain correspondence between weak
limits of powers of a transformation and its skew products.
\end{abstract}

\maketitle

\section*{Introduction}

Let $T$ be a measure-preserving transformation defined on a
non-atomic standard Borel probability space $(X,\mathcal{F} , \mu)$.
The $spectral$ properties of $T$ are those of the induced (
Koopman's) unitary operator on
 $L^2( \mu)$ defined by
\[
\widehat{T}:L^2( \mu)\rightarrow L^2( \mu);\quad
\widehat{T}f(x)=f(Tx).
\]

The transformation $T$ is said to be $\alpha$-rigid on a sequence of
integers $k_i$, where $\alpha \in [0; 1]$, if for every measurable
$A$
\begin{equation}
\lim_i\mu (T^{k_i}A\cap A)\geq \alpha \mu (A).
\end{equation}
 All these $\alpha$'s if we vary the sequence
$k_i$ form a closed subset in $[0; 1]$. The transformation $T$ is
said to be $\alpha$-$rigid$ if $\alpha$ is the right point of this
subset.

By $mixing$ $component$ in the spectrum  of the transformation $T$
we mean a non-zero $\widehat{T}$-invariant subspace, say $H$, of
$L^2( \mu)$ such that $\langle\widehat{T}^kf, g\rangle \rightarrow
0$  as $k\rightarrow \infty$ for every $f,g\in H$.

It is said that a transformation $T$ is a $skew$ $product$ over
$T_0$ if $T$ acts on a product probability Borel space $(X\times
Y,\mathcal{F}\otimes \mathcal{C}, \mu\otimes \nu)$ by
\[
T(x,y)=(T_0x, S_x(y)),
\]
where $(S_x )_{x\in X}$ is a family of transformations of
$(Y,\mathcal{C}, \nu)$ such that the map $(x, y)\mapsto S_xy $ is
measurable. A particular case of it is any $G$-$extension$ over
$T_0$; that is, $(Y,\mathcal{C}, \nu)$ is a group $G$ with the Haar
measure $\nu$, and $S_xy=\phi (x)+y$ for some measurable map
($cocycle$) $\phi: X\rightarrow G$. The restriction of $\widehat{T}$
to an invariant space  $H_0=\{f(x,y)\in L^2:\exists g(x)
[f(x,y)=g(x) \ \mu\otimes \nu- \mbox{a.e.}]\}$ is unitarilly
equivalent to $\widehat{T_0}$. One says that $T$ is $relatively$
$mixing$ if $H_0^\bot$ is its mixing component.

It is well known (see, for example, [4], [8]) that both every
$\alpha$-rigid transformation has no mixing components if $\alpha
>1/2$ and is not mixing if $\alpha >0$. In this short note we prove
that for every $\alpha$, where $\alpha\in [0, 1/2]$, there exists an
$\alpha$-rigid transformation whose spectrum has Lebesgue component.
This paper is based on well-known facts to the author, and we offer
a written version of that because of the current interest to this
theme (see [1]). Recently, the author has been informed about [1],
where, actually, the main result (as it is stated in the abstract)
is that there exists an $\alpha$-rigid transformation with $\alpha$
less or equal to $1/ 2$ whose spectrum has Lebesgue component. Let
me note that the author of [1] knew that there exists an alternative
proof of this theorem. This proof written below is of the
independent interest for possible applications, because, in
particular, of its large constructive potential.

It is well known that the notion to be $\alpha$-rigid transformation
$T$ can be rewritten in terms of the powers weak operator closure
$PCL(T)$ of $\widehat{T}$. Indeed, (0.1) is equivalent to
\[
 \forall f\geq 0  \ \ \lim_i \langle \widehat{T}^{k_i}f, f\rangle\geq \alpha \langle f,
f\rangle.
\]
Taking a subsequence of $k_i$, we can assume that there exist
$Q=\lim_i\widehat{T}^{k_i}$. Obviously, $Q\geq 0$; that is, $Qf\geq
0$ for every $f\geq 0$.

 Let $E$ be the identity map. If for some $f
> 0$ a.e., $Q\prime f= (Q-\alpha \widehat{E} )f < 0 $ on some measurable $A$,
 then
 \[0\leq \langle Q (f\chi_{\overline{A}}),
  f\chi_A\rangle=\langle Q\prime (f\chi_{\overline{A}}),
   f\chi_A\rangle \leq \langle Q\prime f,
   f\chi_A\rangle \leq 0.
    \] Thus $A$ has a zero measure, and $Q\prime\geq 0$. This
    implies that
    \[
    \alpha=\max_{P\in PCL(T)}\{\beta\in [0,1] :\exists Q\geq 0 \ [P= \beta
    \widehat{E}+Q]\}.
    \]

\section{The basic result, applications,  and comments}

It is convenient to think that every $P\in PCL(T_0)$ is also the
operator $P\otimes \widehat{E}$ acting on $L^2(X\times Y,
\mu\otimes\nu)$.

In this paper we prove the following theorem:

\begin{theorem} Let $T$ be a skew product over $T_0$ and $P_{H_0}$
 the orthogonal projection onto $H_0$. If $T$ is relatively mixing, then
\[
PCL(T)=PCL(T_0)P_{H_0}.
\]

\end{theorem}
\begin{proof}
Take $P\in PCL(T)$, $f,g\in L^2$. Denote by $f_0$ and $f_0^\perp$
the corresponding parts of $f$ in $H_0$ and $H_0^\perp$
consequently. Then
\[
\langle Pf, g\rangle= \lim_i\langle \widehat{T}^{k_i}f, g\rangle =.
\lim_i\langle \widehat{T}^{k_i}(f_0+f_0^\perp), g_0+g_0^\perp\rangle
=\lim_i\langle \widehat{T}^{k_i}f_0, g_0\rangle=
\]
\[\lim_i\langle
\widehat{T}_0^{k_i}f_0, g_0\rangle=\lim_i\langle
\widehat{T}_0^{k_i}f_0, g\rangle=\lim_i\langle
\widehat{T}_0^{k_i}P_{H_0}f, g\rangle.\]

Thus $PCL(T)\subseteq PCL(T_0)P_{H_0}$. Along the same line we also
get $PCL(T_0)P_{H_0}\subseteq PCL(T)$.

\end{proof}

\begin{corollary}Let $T$ be a relatively mixing $\textbf{Z}_2$-extension
over $T_0$. Then $T_0$ is $\alpha$-rigid if and only if $T$ is
$\alpha/2$-rigid.
\end{corollary}
\begin{proof} Let us remark first that in this case
$P_{H_0}=(\widehat{E}+\widehat{S})/2$, where $S(x,y)=(x,y+1 \mbox{
mod } 1)$. Applying Theorem next all we have to show is if $\alpha
\widehat{E}+P=1/2Q+1/2Q\widehat{S}$, where $P\geq 0$, $Q\in
PCL(T_0)$, then $Q=2\alpha \widehat{E}+Q\prime$ for some $Q\prime
\geq 0$. Taking $f\geq 0$, $\mbox{supp } f \subseteq \{(x,y):
y=0\}$, ($f\in L^2$), we see that $\langle 1/2Q f, f\rangle=\langle
1/2Q+1/2Q\widehat{S}f, f\rangle \geq \alpha\langle f,f\rangle$. This
implies $\langle 1/2Q f, f\rangle \geq \alpha\langle f,f\rangle$ for
every $f\geq 0$. Thus
 $Q\prime\geq 0$ (see the discussion near definitions of the $\alpha$-rigidity).

   \end{proof}
\begin{corollary}For every $\alpha$, where $\alpha\in [0, 1/2]$,
there exists an $\alpha$-rigid transformation whose spectrum has
Lebesgue component.\footnote{So, it is nonsingular in the sense its
maximal spectral type is nonsingular.} In fact, for every such
$\alpha$ we can take continuum mutually non-isomorphic metrically or
even spectrally $\alpha$-rigid transformations with Lebesgue
component in their spectrum.
\end{corollary}
\begin{proof}Thanks to Helson and Parry (see [6]) every
transformation admits a $\textbf{Z}_2$-extension such that
$H_0^\perp$ is a Lebesgue component in the spectrum. By the
Riemann-Lebesgue lemma every Lebesgue component is mixing. It
remains to mention that for every $\alpha\in [0, 1]$ we can
construct enough of different $\alpha$-rigid transformations
applying the technique developed in [2] or [7].

\end{proof}
\begin{remark} In [3] it is constructed a family of different $\alpha$-rigid
 transformations admitting relatively mixing
$\textbf{Z}_2$-extensions with twofold Lebesque component, where
$\alpha=1$ or less. This implies that there exist 1/2(or less)-rigid
transformations with twofold Lebesque component in its spectrum.
Since Mathew-Nadkarni transformations are included in that family,
and they are $\textbf{Z}_2$-extensions over the $adic$ shift, they
are 1/2-rigid transformations with twofold Lebesque component.
Indeed, every adic shift has the discrete spectrum, therefore it is
just (1-)rigid, and we apply Theorem 1.1.
\end{remark}
\begin{remark} All the below notes admit natural extensions to
both (non)relatively mixing $G$-extensions if we can calculate
powers limits on $H_0^\bot$ and  group actions.
\end{remark}

\bibliographystyle{amsplain}

\end{document}